\documentclass{elsarticle}

\usepackage{lineno,hyperref}

\usepackage{amsmath,amssymb,amsfonts,amsthm}

\newtheorem{theorem}{Theorem}
\newtheorem{remark}{Remark}
\newdefinition{definition}{Definition}
\modulolinenumbers[5]

\journal{Applied Numerical Mathematics}

\newcommand{\RK}{Runge--Kutta}
\newcommand{\RKN}{Runge--Kutta--Nystr\"{o}m}

\begin{document}

\begin{frontmatter}

\title{Functional Continuous Runge--Kutta \\Methods with Reuse}
\author{Alexey S. Eremin}
\address{Department of Information Systems, Saint-Petersburg University, St. Petersburg, 199034, Russia}
\ead{a.eremin@spbu.ru}

\begin{abstract}
In the paper explicit functional continuous Runge--Kutta and Runge--Kutta--Nystr\"om methods for retarded functional differential equations are considered. New methods for first order equations as well as for second order equations of the special form are constructed with the reuse of the last stage of the step. The order conditions for Runge--Kutta--Nystr\"om methods are derived. Methods of orders three, four and five which require less computations than the known methods are presented. Numerical solution of the test problems confirm the convergence order of the new methods and their lower computational cost is performed.
\end{abstract}

\begin{keyword}
functional differential equations
\sep 
continuous Runge--Kutta
\sep
overlapping
\sep
delay differential equations
\sep
\MSC[2010] 65L03\sep 65L06
\end{keyword}

\end{frontmatter}

\section{Introduction}
The paper is in a lot of moments based on the paper by S.~Maset, L.~Torelli and R.~Vermiglio on Functional Continuous Runge--Kutta methods~\cite{MasTorVer2005}. Let's start with some basic denotations.

\begin{itemize}
 \item Let $r \in [0, \infty]$ and $\mathcal{C}$ be the space of continuous functions $[-r,0]\to\mathbb{R}^d$ equipped with the maximum norm
 $$
 \|\varphi\|=\operatorname*{\max}_{\theta\in[-r,\,0]}|\varphi(\theta)|,\quad \varphi\in\mathcal{C},
 $$
where $|\cdot|$ is an arbitrary norm on $\mathbb{R}^d$.

 \item The analogous space of continuously differentiable functions is denoted $\mathcal{C}^1$.

 \item For a continuous function $u:[a-r,b)\to\mathbb{R}^d$ and $t \in [a,b)$, where $a < b$, let $u_t$ be the function given by
\begin{equation}
  u_t(\theta) = u(t + \theta),\quad \theta\in[-r, 0].
  \label{eq:u_t_def}
\end{equation}
\end{itemize}

A differential equation where the higher derivative depends on the unknown function and lower derivatives values in the past is called a \textit{retarded functional differential equation} (RFDE). For example, a first order RFDE is
\begin{equation}
  \dot{u}(t) = f(t, u_t),
  \label{eq:1order}
\end{equation}
and a second order RFDE is
\begin{equation}
  \ddot{u}(t) = f(t, u_t, \dot{u}_t),
\label{eq:2ordergen}
\end{equation}
where 
$$
  \dot{u}_t(\theta) = \dot{u}(t + \theta),\quad \theta\in[-r', 0].  
$$
Various particular cases of RFDEs include \textit{delay differential equations}
$$
  \dot{u}(t) = f(t, u(t), u(t-\tau_1),\dots,u(t-\tau_k)),
  \quad \tau_i \geq 0,\quad i=1,\dots,k,
$$
where in every moment $t$ only the values of $u$ in a finite number of the points in the past are necessary, \textit{integral differential equations}
$$
  \dot{u}(t) = f(t, u(t)) + \int\limits_{t-r}^tK(s,t,u(s))ds,
$$
their combinations or other ways to use the past values $u$.

In the present paper we consider the equation \eqref{eq:1order} and a particular case of \eqref{eq:2ordergen}
\begin{equation}
  \ddot{u}(t) = f(t, u_t),
\label{eq:2order}
\end{equation}

To find a unique solution of an RFDE the initial value is not enough and a \textit{history} function $\varphi$, determining the solution in some interval left of the initial point, is required. In the most cases even for smooth enough $f$ and $\varphi$ the solution doesn't smoothly continue the history. This leads to a number of points where the solution has jump discontinuities in some derivatives. This restricts greatly multistep methods application and in recently the main attention was devoted to one-step methods, specifically continuous Runge--Kutta (CRKs)~\cite{BellZenn}.

A CRK provides a continuous approximation of the solution over the integration step, which can be later substituted into the right-hand side $f$ when needed. However, in the case, when we need the continuous approximation within the currently calculated step, the implementation of any (even explicit) Runge--Kutta method becomes fully implicit. This situation is known as \textit{overlapping} when delay differential equations are considered. For integral differential equations or more general types of RFDEs such situation occurs at every step.

Overlapping makes application of Bellmann's method of steps~\cite{BelCoo1965} or explicit Runge--Kutta methods impossible. Though fully implicit methods (like RADAR code by Guglielmi and Hairer~\cite{GugHai2008}) work fine, still the speed of explicit methods is often desirable.

A way to construct explicit methods for general RFDEs was first proposed by Tavernini in early seventies~\cite{Tavernini1971} but only few decades later his approach was further developed by a group of Italian researchers~\cite{MasTorVer2005, ActaNumer2009}. They provide the continuous approximations of rising orders for every stage, finally reaching the desired method's order. Such methods, named Functional Continuous Runge--Kutta methods (FCRKs), are the subject of the present paper. 

It should be mentioned, that general linear multistep methods are studied as a way to solve RFDEs as well (e.g.,~\cite{Pimenov1999}), and even
a functional continuous approach is used in them~\cite{Tuzov2017Preprint}. Still due to the reasons mentioned above we concentrate on one-step methods.

The methods constructed in~\cite{MasTorVer2005} can be made less expensive if one uses the last stage of the step as the first stage of the next step, as it was done for instance for CRKs in~\cite{OwrZen1992}. 

In the next section we recall the necessary information on FCRKs, and then in Sec.~\ref{sec:fcrk_reuse} construct FCRKs with the last stage reuse. We also study FCRK methods for direct application to the second order equations of special form, which are analogous to Runge--Kutta--Nystr\"om methods (Sec.~\ref{sec:fcrkn}), prove their order conditions (Sec.~\ref{sec:fcrkn_ocs}) and finally present such methods with reuse (Sec.~\ref{sec:fcrkn_reuse}). In the last section we run test problems that demonstrate the convergence of the presented methods.

\section{Runge--Kutta Methods for RFDEs \label{sec:fcrk}} 
This section recalls the results presented in~\cite{MasTorVer2005}.
We consider only explicit method in the current paper and make the corresponding changes to the cited material. 

Here and in the next section we consider the first order RFDE 
\begin{equation*}
  \dot{u}(t) = f(t, u_t),
  \tag{\ref{eq:1order}}
\end{equation*}
where $f: \Omega\to\mathbb{R}^d$, and open set $\Omega\subseteq\mathbb{R}\times\mathcal{C}$. We assume that $f$ is continuous and its derivative $f' : \Omega \to \mathcal{L}(\mathcal{C}, \mathbb{R}^d)$ is bounded and continuous with respect to the second argument. In this case according to \cite{HaleVerLunel} for each $(\sigma, \varphi) \in \Omega$ there exists a unique (non-continuable) solution $u = u(\sigma, \varphi) : [\sigma - r, \bar{t}) \to \mathbb{R}^d$ of \eqref{eq:1order} through $(\sigma, \varphi)$, where $\bar{t} = \bar{t}(\sigma, \varphi) \in (\sigma, +\infty]$, i.e. $u$ satisfies \eqref{eq:1order} for $t\in [\sigma, \bar{t})$ and $u_\sigma = \varphi$.

\begin{definition}
Let $s$ be a positive integer. An explicit $s$-stage functional continuous Runge--Kutta method (FCRK) is a triple ($A(\cdot)$, $b(\cdot)$, $c$) where
\begin{itemize}
 \item $A(\cdot)$ is a strict lower-triangular $\mathbb{R}^{s\times s}$-valued polynomial function such that $A(0) = 0$,
 
 \item $b(\cdot)$ is an $\mathbb{R}^s$-valued polynomial function such that $b(0)=0$,
 
 \item $c \in \mathbb{R}^s$ with $c_1 = 0$ and $c_i \geq 0$, $i = 2,..., s$.
\end{itemize}
\end{definition}

Applied with stepsize $h$ to \eqref{eq:1order} to get the solution $u$ through $(\sigma,\varphi)$, the FCRK ($A(\cdot)$, $b(\cdot)$, $c$) provides the continuous approximation $\eta(\alpha h)$ of the shift $y=u(\sigma+\cdot)$ on $[0, h]$:
\begin{equation}
  \eta(\alpha h) = \varphi(0) + h\sum_{i=1}^sb_i(\alpha)K_i,\quad \alpha\in[0, 1],
  \label{eq:fcrk_u}
\end{equation}
where
\begin{equation}
  K_i = f\!\left(\sigma + c_ih, Y_{c_ih}^{i}\right)\!,\quad i=1,..., s
  \label{eq:fcrk_k}
\end{equation}
and $Y^{i} : [-r, c_ih] \to \mathbb{R}^d$ are \textit{stage functions} given by 
\begin{equation}
\begin{aligned}
  Y^{i}(\alpha h) &= \varphi(0) + h\sum_{j=1}^{i-1}a_{ij}(\alpha)K_j, &&\alpha\in[0,c_i],\\
  Y^{i}(\theta) &= \varphi(\theta),&& \theta\in[-r,0].
\end{aligned}
\label{eq:fcrk_Y}
\end{equation}
The conditions $A(0) = 0$ and $b(0) = 0$ guarantee $Y_{c_ih}^{i}\in\mathcal{C}$, $i = 1,...,s$,
and $\eta_h\in \mathcal{C}$ respectively.

When the second step is made, the function $\varphi$ in \eqref{eq:fcrk_u}--\eqref{eq:fcrk_Y} is extended up to the new starting point ($\sigma+h$) with $\eta$ from the first step. The same for the following steps.

\begin{definition}
The function
\begin{equation}
  E = E(h,\sigma,\varphi) = \eta - y: \quad [0, h] \to \mathbb{R}^d,
  \label{eq:locerr}
\end{equation}
is called the \textit{local error}. We say that for a sufficiently smooth problem an FCRK has \textit{local uniform (discrete) order} $p$ ($q$) if for $h$ small enough there exists some $C>0$ such that
\begin{equation*}
  \max_{\alpha\in[0,1]}\|E(\alpha h)\|\leq Ch^{p+1}\quad
  \left(\|E(h)\|\leq Ch^{q+1}\right).
\end{equation*}
\end{definition}
It is obvious, that $q\geq p$. More rigorous definitions, which take in account discontinuity points, can be found in~\cite{MasTorVer2005}. The problem of global convergence and its connection to the local orders is considered in~\cite{BellZenn}. It is enough to mention here that a method needs to have discrete order $p$ and uniform order $p-1$ to provide the convergence order $p$. Still we construct uniform order $p$ methods here, since when implemented they are better in various senses (more justified local error estimation, its minimization, application to neutral equations, etc.)

An FCRKs can be conveniently presented with a Butcher tableau:
\begin{equation}
  \begin{array}{c|ccccc}
    0    & & &  \\
    c_2  & a_{21}(\alpha) \\
    c_3  & a_{31}(\alpha) & a_{32}(\alpha) & \\
    \vdots & \vdots & \vdots & \ddots  & \\
    c_s  & a_{s1}(\alpha) & a_{s2}(\alpha) & \cdots & a_{s,s-1}(\alpha) \\\hline
         & b_1(\alpha) & b_2(\alpha) & \cdots & b_{s-1}(\alpha) & b_s(\alpha)
  \end{array}
  \label{tbl.frck_gen}
\end{equation}
It can be reduced to a continuous Runge--Kutta (CRK) method for ODEs (or DDEs with non-vanishing delays) by setting $a_{ij}=a_{ij}(c_i)$.

In~\cite{MasTorVer2005} the methods of uniform orders 1, 2, 3 and 4 were presented with 1, 2, 4 and 7 stages respectively. Those are the lowest numbers of stages providing such \textit{uniform} orders. However, it is possible to construct a \textit{discrete} order 3 method with 3 stages and a \textit{discrete} order 4 method with 6 stages. This leads to methods with \textit{reuse} studied in the next section.

\section{Methods with the Last Stage Reuse \label{sec:fcrk_reuse}}
Continuous \RK{} methods (CRKs) are extensions of Runge--Kutta methods for an ODE initial value problem 
$$
  \dot{u}(t) = f(t, u(t)), \quad y(t_0)=y_0
$$
providing the continuous approximation $\eta(\alpha h)$ of the solution $u(t)$ on $[t_0,t_0+h]$
\begin{equation}
\begin{aligned}
 \eta(\alpha h) &= y_0 + h\sum_{i=1}^sb_i(\alpha)K_i,\quad \alpha\in[0, 1],\\
 K_i &= f\!\left(\sigma + c_ih, Y_i\right)\!, \quad Y_i = y_0 + h\sum_{j=1}^{i-1}a_{ij}K_j, \quad i=1,..., s.
\end{aligned}
\end{equation}
Since here the matrix $A$ is constant those methods have less strict order conditions than FCRKs, and thus for orders 4 and higher can be constructed with fewer stages~(see \cite{BellZenn}). However, since they find wide application in solution of DDEs (and also RFDEs) with non-vanishing delays when one uses smaller step sizes than the minimum delay value, they are usually constructed to have the uniform order equal to the discrete order, or at least 1 order lower. Owren and Zennaro~\cite{OwrZen1992} have constructed ``optimal'' CRKs of orders up to five in which the idea of getting the discrete order $p$ with fewer stages than it is necessary for uniform order $p$ is used. The additional stage necessary for uniform order $p$ is computed in the point $(t_0+h, \eta(h))$ (which is order $p$ approximation to $u(t_0+h)$) and can be used as a first stage for the next step~--- the approach named \textit{reuse} or \textit{First Same as Last, FSAL}.
We don't recall details on CRKs here. They can be found in the cited works~\cite{BellZenn, OwrZen1992}. Let's show how the same idea can be applied for FCRKs.

As it was already mentioned, methods of discrete orders 3 and 4 can be constructed with just 3 and 6 stages. In both cases one additional stage is sufficient to provide the uniform order 3 or 4 as well. This last stage will be used a first stage of the next step. 

The general formulation of the method remains the same as \eqref{eq:fcrk_u}--\eqref{eq:fcrk_Y}. We only have additional restrictions on the parameters:
\begin{itemize}    
  \item $c_s=1$;
  \item $b_i(1) = a_{si}(1)$ for any $i=1,...,s$;
  \item $a_{si}(\alpha)$ must satisfy discrete order $p$ and uniform order $p-1$ conditions as $b$-parameters of a method with $s-1$ stages.
\end{itemize}
The first condition is necessary to reuse the stage, while the second one provides that the continuous extension $\eta$ ends in the point obtained by the method of discrete order $p$ with $s-1$ stages. 

We must not only provide the discrete order $p$ with $s-1$ stages, but the uniform order $p-1$ as well. This is necessary due to the fact, that the low order of the last stage of the previous step can reduce the order at the current step.

It should also be noted that if the step starts from the point where $\dot{u}$ has a jump discontinuity (which for DDEs can occur only for the first step, or if the history $\varphi$ or right-hand side $f$ have jumps), we don't use the last stage of the previous step and recompute it with the new branch of $\varphi$ or $f$. Some details on discontinuity approximation and branch-wise control of problems smoothness can be found in~\cite{EreHum2015}.

We now present methods of orders 3 and 4.

\subsection{Method of order three}
The method obtains the solution in the next mesh point with 3 stages and uses the value $K_4$ to get order 3 uniform approximation. Free parameters are chosen to reduce the error coefficients of order four in mean square sense (they were computed only for application of the method to ODEs or DDEs with nonvanishing delays).
\begin{equation}
  \begin{array}{c|cccc}
    0    & & &  \\
    \frac12    & \alpha \\
    \frac23    & \alpha-\alpha^2 & \alpha^2 \\
    1    & \alpha-\frac34\alpha^2 & 0  & \frac34\alpha^2 \\\hline
      & \alpha-\frac54\alpha^2+\frac12\alpha^3 & 0 &  \frac94\alpha^2-\frac32\alpha^3 & -\alpha^2+\alpha^3
  \end{array}
  \label{tbl.fcrk3f}
\end{equation}

\subsection{Method of order four}
Here only six new stages are required for every step. Free parameters were chosen to reduce the error as well as for the method of order 3.
\begin{equation}
  \begin{array}{c|ccccccc}
    0    & & &  \\%
    \frac25       & \alpha \\%
    \frac7{19}    & \alpha-\frac54\alpha^2 & \frac54\alpha^2 \\%
    \frac{15}{17} & \alpha-\frac54\alpha^2 & \frac54\alpha^2 \\%
    \frac{5}{14}  & a_{51}(\alpha) & 0 &  a_{53}(\alpha) & a_{54}(\alpha) \\%
    \frac{11}{13} & a_{61}(\alpha) & 0 &  a_{63}(\alpha) & a_{64}(\alpha)   \\%
    1 & a_{71}(\alpha) & 0 &  0 & 0 & a_{75}(\alpha) & a_{76}(\alpha) 
    \\\hline%
      & b_1(\alpha) & 0 & 0 & 0 & b_5(\alpha) & b_6(\alpha) & b_7(\alpha)
  \end{array}
  \label{tbl.fcrk4f}
\end{equation}
$$
\begin{aligned}
  &a_{51}(\alpha) = a_{61}(\alpha)=\alpha-\frac{202}{105}\alpha^2+\frac{323}{315}\alpha^3, \\
  &a_{53}(\alpha) = a_{63}(\alpha) = \frac{5415}{2324}\alpha^2-\frac{6137}{3486}\alpha^3, \\
  &a_{54}(\alpha) = a_{64}(\alpha) = -\frac{2023}{4980}\alpha^2+\frac{5491}{7470}\alpha^3,
  \end{aligned}
$$
$$
\begin{aligned}
  &a_{71}(\alpha) = \alpha-\frac{219}{210}\alpha^2+\frac{182}{165}\alpha^3,  
   && b_1(\alpha) = \alpha-\frac{137}{55}\alpha^2+\frac{401}{165}\alpha^3-\frac{91}{110}\alpha^4,
  \\
  &a_{75}(\alpha) = \frac{1078}{445}\alpha^2-\frac{2548}{1335}\alpha^3,
  && b_5(\alpha) = \frac{15092}{4005}\alpha^2
  -\frac{21952}{4005}\alpha^3 +\frac{8918}{4005}\alpha^4,
  \\
  &a_{76}(\alpha) = -\frac{845}{1958}\alpha^2+\frac{2366}{2937}\alpha^3,
  && b_6(\alpha) = -\frac{10985}{3916}\alpha^2+\frac{41743}{5874}\alpha^3-\frac{15379}{3916}\alpha^4,
  \\
  &&&b_7(\alpha) = \frac{55}{36}\alpha^2-\frac{73}{18}\alpha^3+\frac{91}{36}\alpha^4.
\end{aligned}
$$

\section{\RKN{} Methods for Second Order Equations \label{sec:fcrkn}}
Here we consider the equation
\begin{equation}
 \ddot{u}(t) = f(t,u_t).
 \tag{\ref{eq:2order}}
\end{equation}
The solution existence and uniqueness conditions for it can be obtained by rewriting it as a first order system and applying results from~\cite{HaleVerLunel} as for \eqref{eq:1order}. Namely, 
we assume that $u_t\in\mathcal{C}^1$, $f : \Omega\to\mathbb{R}^d$ and $\Omega$ is an open subset of $\mathbb{R}\times\mathcal{C}^1$, 
$f$ is continuous and has derivative $f' : \Omega \to \mathcal{L}(\mathcal{C}^1, \mathbb{R}^d)$ with respect to the second argument which is bounded and continuous with respect to the second argument. 
Thus, for each $(\sigma, \varphi) \in \Omega$ there exists a unique (non-continuable) solution $u = u(\sigma, \varphi) : [\sigma - r, \bar{t}) \to \mathbb{R}^d$ of \eqref{eq:2order} through $(\sigma, \varphi)$, where $\bar{t} = \bar{t}(\sigma, \varphi) \in (\sigma, +\infty]$, i.e.\ $u$ satisfies \eqref{eq:2order} for $t\in [\sigma, \bar{t})$ and $u_\sigma = \varphi$.

Notice that since the right-hand side of \eqref{eq:2order} doesn't depend on $\dot{u}$ there is no need to need in any additional assumptions on $\dot{\varphi}$, save the existence of the initial value $\dot{\varphi}(0)$.

\begin{remark} Through the whole paper we use dot (\,$\dot{}$\,) for time (or time-like variable) derivative, and upper index in brackets ($^{(k)}$) means $k$-th time derivative. Prime (\,$'$\,) is only used  to indicate that the function or parameter is somehow connected to dot-variables and is specific for second order equations, i.e. prime does \textit{not} mean derivation. 
\end{remark}

In full analogy with FCRK \eqref{eq:fcrk_u}--\eqref{eq:fcrk_Y} the following method for direct implementation to \eqref{eq:2order} we introduced in~\cite{EJQTDE2015}.

\begin{definition}
Let $s$ be a positive integer. An explicit $s$-stage functional continuous Runge--Kutta--Nystr\"om method (FCRKN) is a quadruple ($A(\cdot)$, $b(\cdot)$, $b'(\cdot)$, $c$) where
\begin{itemize}
 \item $A'(\cdot)$ is a strict lower-triangular $\mathbb{R}^{s\times s}$-valued polynomial function such that $A(0) = 0$,

 \item $b(\cdot)$ and $b'(\cdot)$ are $\mathbb{R}^s$-valued polynomial functions such that $b(0)=b'(0)=0$,

 \item $c \in \mathbb{R}^s$ with $c_1 = 0$ and $c_i \geq 0$, $i = 2,..., s$.
\end{itemize}
\end{definition}

Applied with stepsize $h$ to \eqref{eq:2order} to get the solution $u$ through $(\sigma,\varphi)$, the FCRKN ($A(\cdot)$, $b(\cdot)$, $b'(\cdot)$, $c$) provides the continuous approximations $\eta(\alpha h)$ of the shift $y=u(\sigma+\cdot)$ and $\eta'(\alpha h)$ of the shift $\dot{y}=\dot{u}(\sigma+\cdot)$ on $[0, h]$:
\begin{equation}
  \begin{aligned}
    \eta(\alpha h) &= \varphi(0) + \alpha h \dot{\varphi}(0) + h^2\sum_{i=1}^sb_i(\alpha)K_i, &&\alpha\in[0, 1],\\
    \eta'(\alpha h) &= \dot{\varphi}(0) + h\sum_{i=1}^sb'_i(\alpha)K_i,
  \end{aligned}
  \label{eq:fcrkn_u}
\end{equation}
where again
\begin{equation}
  K_i = f\!\left(\sigma + c_ih, Y_{c_ih}^{i}\right)\!,\quad i=1,..., s
  \label{eq:fcrkn_k}
\end{equation}
and $Y^{i} : [-r, c_ih] \to \mathbb{R}^d$ are \textit{stage functions} given by 
\begin{equation}
\begin{aligned}
  Y^{i}(\alpha h) &= \varphi(0) + \alpha h \dot{\varphi}(0) + h^2\sum_{j=1}^{i-1}a_{ij}(\alpha)K_j, &&\alpha\in[0,c_i],\\
  Y^{i}(\theta) &= \varphi(\theta),&& \theta\in[-r,0].
\end{aligned}
\label{eq:fcrkn_Y}
\end{equation}
The conditions $A(0) = 0$, $b(0) = 0$ and $b'(0)$ guarantee $Y_{c_ih}^{i}\in\mathcal{C}$, $i = 1,...,s$,
$\eta_h\in \mathcal{C}$ and $\eta'_h\in \mathcal{C}$ respectively.

\begin{remark} The coefficients $b'_i(\alpha)$ of the FCRKN \eqref{eq:fcrkn_u} are connected to the $b_i(\alpha)$ coefficients of the FCRK \eqref{eq:fcrk_u} closer than $b_i(\alpha)$ of \eqref{eq:fcrkn_u} are. In fact an application of FCRK to the system
$$
  \left\{\begin{aligned}
   \dot{u}(t) &= v(t),\\
   \dot{v}(t) &= f(t,u_t),
   \end{aligned}\right.
$$
which is equivalent to \eqref{eq:2order} will lead to the approximation of $\dot{u}$ in the form
$$
  \begin{aligned}
    v(\sigma+\alpha h) \approx \eta'(\alpha h) &= \dot{\varphi}(0) + h\sum_{i=1}^sb_i(\alpha)K_i.
  \end{aligned}
$$
Still to make the denotations for FCRKNs more consistent, we correspond $b'_i(\alpha)$ to $\eta'(\alpha h)$ and $b_i(\alpha)$ to $\eta(\alpha h)$.

\end{remark}

FCRKNs orders are defined in more complicated way than those of FCRKs. Along with the local error \eqref{eq:locerr} the local error is introduced for $\eta'$:
\begin{equation}
  E' = E'(h,\sigma,\varphi) = \eta' - \dot{y}: \quad [0, h] \to \mathbb{R}^d,
  \label{eq:locerrp}
\end{equation}

\begin{definition}
We say that for a sufficiently smooth problem \eqref{eq:2order} an FCRKN has \textit{local uniform (discrete) order} $p$ ($q$) if for $h$ small enough there exist some $C>0$ and $C'>0$ such that
\begin{gather*}
  \max_{\alpha\in[0,1]}\|E(\alpha h)\|\leq Ch^{p+1} \text{ and } 
  \max_{\alpha\in[0,1]}\|E'(\alpha h)\|\leq C'h^{p+1}\\
  \bigg(\|E(h)\|\leq Ch^{q+1} \text{ and } \|E'(h)\|\leq C'h^{q+1}\bigg).
\end{gather*}
\end{definition}

\section{Order Conditions \label{sec:fcrkn_ocs}}

In~\cite{EJQTDE2015} FCRKNs order conditions were presented without demonstration of their necessity or even sufficiency. Intuitively constructed labeled trees correspondence to order conditions was used. However, as it was mentioned there, a strict proof was still needed. Since a separate paper on FCRKNs order conditions derivation would be almost useless, we include the rigorous order conditions derivation here.

\subsection{Error expansions}
Analogously to \textit{local errors}
\begin{equation*}
 \begin{aligned}
  E &= \eta - y: &&[0, h] \to R^d,\\
  E' &= \eta' - \dot{y}: &&[0, h]\to R^d\\
 \end{aligned}
\end{equation*}
we introduce \textit{stage errors}
\begin{equation}
  E^i = Y^i - y: [-r, c_ih] \to R^d, \quad i=1,...,s.
\end{equation}

We also extend $a_{ij}(\alpha) =0$, $\alpha\leq0$ for all $i,j=1,...,s$ and denote $\bar{a}_{ij}(\omega)=a_{ij}(c_i+\omega)$, $\omega \leq0$.

Let's study local errors
$$
  E(\alpha h) = \eta(\alpha h) - y(\alpha h)
$$
and
$$
  E'(\alpha h) = \eta'(\alpha h) - \dot{y}(\alpha h).
$$

Notice that $Y^i_{c_ih} = y_{c_ih} + E^i_{c_ih}$ and $\ddot{y}(c_ih) = f(\sigma+c_ih, y_{c_ih})$. We also introduce
\begin{equation}
 D_i = f(\sigma+c_ih, y_{c_ih} + E^i_{c_ih})- f(\sigma+c_ih, y_{c_ih}),\quad i=1,...,s.
 \label{D_i}
\end{equation}

First for Nystr\"om methods with $\alpha\in[0,1]$
$$
\begin{aligned}
  E(\alpha h) &= h^2\sum_{i=1}^s b_i(\alpha)f(\sigma+c_ih, Y^i_{c_ih}) + \varphi(0) + h\dot\varphi(0) - y(\alpha h)\\
  &= h^2\sum_{i=1}^s b_i(\alpha)f(\sigma+c_ih, y_{c_ih} + E^i_{c_ih}) + \varphi(0) + h\dot\varphi(0) - y(\alpha h)\\
  &= h^2\sum_{i=1}^s b_i(\alpha)\left[f(\sigma+c_ih, y_{c_ih} + E^i_{c_ih})- f(\sigma+c_ih, y_{c_ih})\right]\\&\hspace{5mm} + h^2\sum_{i=1}^s b_i(\alpha)f(\sigma+c_ih, y_{c_ih}) +  \varphi(0) + h\dot\varphi(0) - y(\alpha h)\\
  &= h^2\sum_{i=1}^s b_i(\alpha)D_i + h^2\sum_{i=1}^s b_i(\alpha)\ddot{y}(c_ih) +  \varphi(0) + h\dot\varphi(0) - y(\alpha h).
\end{aligned}
$$
Now by Taylor expansion of $y$ and $\ddot{y}$ we get
\begin{equation}
  E(\alpha h)= h^2\sum_{i=1}^s b_i(\alpha)D_i + \sum_{k=2}^p h^k \Gamma_k(\alpha)y^{(k)}(0) + O(h^{p+1}),\quad\alpha\in[0,1],
  \label{E_taylor}
\end{equation}
where
\begin{equation}
 \Gamma_k(\alpha) = \frac{1}{(k-2)!}\left(\sum_{i=1}^s b_i(\alpha)c_i^{k-2} - \frac{\alpha^k}{k(k-1)}\right)\!\!,\quad\alpha\in[0,1].
 \label{gamma_k}
\end{equation}

Analogously 
\begin{equation}
  E'(\alpha h)= h\sum_{i=1}^s b'_i(\alpha)D_i + \sum_{k=1}^p h^k \Gamma'_k(\alpha)y^{(k+1)}(0) + O(h^{p+1}),\quad\alpha\in[0,1],
  \label{Ep_taylor}
\end{equation}
where
\begin{equation}
 \Gamma'_k(\alpha) = \frac{1}{(k-1)!}\left(\sum_{i=1}^s b'_i(\alpha)c_i^{k-1} - \frac{\alpha^k}{k}\right)\!\!,\quad\alpha\in[0,1].
 \label{gamma_kp}
\end{equation}

Moreover for stage errors
\begin{equation}
  E^i_{c_ih} = h^2\sum_{j=1}^s \bar{a}_{ij}\left(\frac{\cdot}{h}\right)D_j + \sum_{k=2}^p h^k \bar\Gamma_{ik}\left(\frac{\cdot}{h}\right)y^{(k)}(0) + O(h^{p+1}),\quad i=1,...,s,
  \label{Ei_taylor}
\end{equation}
where $\bar{\Gamma}_{ik}(\omega) = \Gamma_{ik}(c_i+\omega)$, $\omega\leq0$ with 
\begin{equation}
 \Gamma_{ik}(\alpha) = \frac{1}{(k-2)!}\left(\sum_{j=1}^s a_{ij}(\alpha)c_j^{k-2} - \frac{\alpha^k}{k(k-1)}\right)\!\!,\quad\alpha\in[0,c_i], \quad i=1,...,s
 \label{gamma_ik}
\end{equation}
and $\Gamma_{ik}(\alpha) = 0$ for $\alpha\leq0$.

By considering RFDEs with $f(t, u_t) = g(t)$ (pure quadrature problems), it is
easy to see that
$$
 \begin{aligned}
  &\Gamma_k = 0,\;\; k = 2,...,p, &&   \Gamma'_k = 0,\;\; k = 1,...,p,\\  
  \Big(&\Gamma_k(1) = 0,\;\;k = 2,...,p, && \Gamma'_k(1) = 0,\;\;k = 1,...,p\Big)  
 \end{aligned}
$$
are necessary conditions for the uniform (discrete) order $p$.

We assume in the following that FCRKN methods satisfy
\begin{equation}
 \sum_{i=1}^sb'_i(\alpha) = \alpha,\quad \alpha\in[0,1],
 \label{OCN1}
\end{equation}
i.e. $\Gamma'_1=0$, which is the uniform order one condition, and also
\begin{equation}
 \begin{aligned}
  \sum_{i=1}^sb_i(\alpha) = \frac{\alpha^2}{2},&&& \alpha\in[0,1],\\
  \sum_{j=1}^sa_{ij}(\alpha) = \frac{\alpha^2}{2},&&& \alpha\in[0,c_i],&& i=1,...,s, 
  \label{OCN2p}
 \end{aligned}
\end{equation}
i.e. $\Gamma_2=0$ and $\Gamma_{i2}=0$, $i=1,...,s$. The first equation of \eqref{OCN2p} is the uniform order two condition for $y$ approximation and the other are simplifying conditions providing the uniform order two of stage approximations $Y^i$.

\subsection{Second order}
We assumed for the existence of the solution of \eqref{eq:2order} that $f$ is of class $C^1$ with respect to the second argument. Now let us assume also that $u$ is of piecewise class $C^3$.

Under these assumptions it is clear that for a method of order one
\begin{equation}
 D_i = O(h^2),\quad i=1,...,s,
 \label{di_oh2}
\end{equation}
in \eqref{D_i}, and thus
\begin{equation}
 \begin{aligned}
  E(\alpha h) &= O(h^3), && \alpha\in[0,1],\\
  E'(\alpha h) &= h^2\Gamma'_2(\alpha)y^{(3)}(0) + O(h^3), && \alpha\in[0,1].
  \label{e_oh3}
 \end{aligned}
\end{equation}

\begin{theorem}\label{ThOCN2}
  An FCRKN method satisfying \eqref{OCN1} and \eqref{OCN2p} is of uniform order two iff
  \begin{equation}
   \sum_{i=1}^sb'_i(\alpha)c_i = \frac{\alpha^2}{2},\quad\alpha\in[0,1]
   \label{OCN2}
  \end{equation}
  and is of discrete order two iff
  \begin{equation}
   \sum_{i=1}^sb'_ic_i = \frac{1}{2}.
   \label{OCN2d}
  \end{equation}
\end{theorem}
\noindent\textbf{Proof.} Let us observe that \eqref{OCN2} is equivalent to $\Gamma'_2 = 0$ (\eqref{OCN2d} is equivalent to $\Gamma'_2(1) = 0$). The ``if'' part follows by \eqref{e_oh3}. As for the ``only if'' part we have
remarked above that $\Gamma'_2 = 0$ ($\Gamma'_2(1) = 0$) is a necessary condition for uniform (discrete) order two.
\qed

\subsection{Third order}
Now we develop the conditions for uniform and discrete orders three. Let us assume that $u$ is of piecewise class $C^4$.

Because of \eqref{di_oh2} and \eqref{OCN2p}
$$
 E^i_{c_ih} = O(h^3),\quad i=1,...,s
$$
and then
\begin{equation}
 D_i = O(h^3), \quad i=1,...,s.
 \label{di_oh3}
\end{equation}
Thus
\begin{equation}
 \begin{aligned}
  E(\alpha h) &= h^3\Gamma_3(\alpha)y^{(3)}(0) +  O(h^4), && \alpha\in[0,1],\\
  E'(\alpha h) &= h^2\Gamma'_2(\alpha)y^{(3)}(0) + h^3\Gamma'_3(\alpha)y^{(4)}(0) + O(h^4), && \alpha\in[0,1].
 \end{aligned}
 \label{e_oh4}
\end{equation}

\begin{theorem}
  An FCRKN method satisfying \eqref{OCN1} and \eqref{OCN2p} and  of uniform order two is of uniform order three iff
  \begin{equation}
   \sum_{i=1}^sb'_i(\alpha)c_i^2 = \frac{\alpha^3}{3},\quad\alpha\in[0,1]
   \label{OCN3}
  \end{equation}
  and 
  \begin{equation}
   \sum_{i=1}^sb_i(\alpha)c_i = \frac{\alpha^3}{6},\quad\alpha\in[0,1],
   \label{OCN3p}
  \end{equation}
  and is of discrete order three iff
  \begin{equation}
   \sum_{i=1}^sb'_ic_i^2 = \frac{1}{3}
   \label{OCN3d}
  \end{equation}
  and
  \begin{equation}
   \sum_{i=1}^sb_ic_i = \frac{1}{6}.
   \label{OCN3pd}
  \end{equation}
\end{theorem}
\noindent\textbf{Proof.} The proof is as straighforward as of the Th.\,\ref{ThOCN2}. Observe that \eqref{OCN3} is equivalent to $\Gamma'_3 = 0$ (\eqref{OCN3d} is equivalent to $\Gamma'_3(1) = 0$) and \eqref{OCN3p} is equivalent to $\Gamma_3 = 0$ (\eqref{OCN3pd} is equivalent to $\Gamma_3(1) = 0$). The ``if'' part follows by \eqref{e_oh4}. Again the ``only if'' part is provided by the fact that $\Gamma'_3 = 0$ ($\Gamma'_3(1) = 0$) and $\Gamma_3 = 0$ ($\Gamma_3(1) = 0$) are necessary conditions for uniform (discrete) order three.
\qed

\subsection{Fourth order}
Now let us assume that $f$ is of class $C^2$ with respect to second argument and $u$ is of piecewise class $C^5$.

Because of \eqref{di_oh3} and \eqref{OCN2p}
$$
 E^i_{c_ih} = h^3\bar{\Gamma}_{i3}\left(\frac{\cdot}{h}\right)y^{(3)}(0) + O(h^4),\quad i=1,...,s
$$
and then
\begin{equation}
 D_i = h^3\left.\frac{\delta f}{\delta u}(\sigma+c_ih,y_{c_ih})\right|_{\theta}\left(\bar{\Gamma}_{i3}\left(\frac{\theta}{h}\right)y^{(3)}(0)\right) + O(h^4), \quad i=1,...,s,
 \label{di_oh4}
\end{equation}
where the symbol $|_{\theta}$ shows that the functional derivative~\cite{GelfFomin} $\frac{\delta f}{\delta u}(\sigma+c_ih,y_{c_ih})$  is applied to the function $\bar{\Gamma}_{i3}\left(\frac{\theta}{h}\right)y^{(3)}(0)$ of the variable $\theta$, $\theta\in[-r,0]$.

Now we have 
\begin{equation}
  E(\alpha h) = h^3\Gamma_3(\alpha)y^{(3)}(0) +  h^4\Gamma_4(\alpha)y^{(4)}(0) +  O(h^5), \quad \alpha\in[0,1],
  \label{e_oh5}
\end{equation}  
\begin{equation}
 \begin{aligned}
  E'(\alpha h) &= h^2\Gamma'_2(\alpha)y^{(3)}(0) + h^3\Gamma'_3(\alpha)y^{(4)}(0) \\
  &\hspace{5mm}+ h^4 \sum_{m=1}^{s^*}\left.\frac{\delta f}{\delta u}(\sigma+c^*_mh,y_{c^*_mh})\right|_{\theta}\left(\sum_{\substack{i=1\\c_i=c^*_m}}^sb'_i(\alpha)\bar{\Gamma}_{i3}\left(\frac{\theta}{h}\right)y^{(3)}(0)\right)\\
  &\hspace{5mm}+ h^4\Gamma'_4(\alpha)y^{(5)}(0) + O(h^5), \quad \alpha\in[0,1],
 \end{aligned}
 \label{ep_oh5}
\end{equation}
where the sum in brackets is made only for $i$ for which $c_i = c^*_m$.

\begin{theorem}\label{TH_OCN4}
  An FCRKN method satisfying \eqref{OCN1} and \eqref{OCN2p} and  of uniform order three is of uniform order four iff
  \begin{equation}
   \sum_{i=1}^sb'_i(\alpha)c_i^3 = \frac{\alpha^4}{4},\quad\alpha\in[0,1],
   \label{OCN41}
  \end{equation}
  \begin{equation}
   \sum_{\substack{i=1\\c_i=c^*_m}}^sb'_i(\alpha)\left(\sum_{j=1}^{s}a_{ij}(\beta)c_j - \frac{\beta^3}{6}\right)\!\!,\quad\alpha\in[0,1],\quad \beta\in[0,c^*_m]
   \label{OCN42}
  \end{equation}
  for $m=1,...,s^*$, and 
  \begin{equation}
   \sum_{i=1}^sb_i(\alpha)c_i^2 = \frac{\alpha^4}{12},\quad\alpha\in[0,1].
   \label{OCN41p}
  \end{equation}
\end{theorem}
\noindent\textbf{Proof.} The proof is analogous to the order three conditions for FCRKs for a first-order RFDE given in~\cite{MasTorVer2005}.

For $m = 1,...,s^*$, let $\Delta_{m3}$ be the function given by
$$
\Delta_{m3}(\alpha, \omega) =
\sum_{\substack{i=1\\c_i=c^*_m}}^sb'_i(\alpha)\bar{\Gamma}_{i3}(\omega),\quad \alpha\in[0, 1], \quad\omega \leq 0.
$$
Then \eqref{ep_oh5} can be written as
$$
\begin{aligned}
  E'(\alpha h) &= h^4 \sum_{m=1}^{s^*}\left.\frac{\delta f}{\delta u}(\sigma+c^*_mh,y_{c^*_mh})\right|_{\theta}\left(\Delta_{m3}\left(\alpha, \frac{\theta}{h}\right)y^{(3)}(0)\right)\\
  &\hspace{5mm}+ h^4\Gamma'_4(\alpha)y^{(5)}(0) + O(h^5), \quad \alpha\in[0,1],
 \end{aligned}
$$
under the assumption that the method has uniform order three.
Let us observe that \eqref{OCN41} is equivalent to $\Gamma'_4 = 0$, \eqref{OCN41p} is equivalent to $\Gamma_4 = 0$ and \eqref{OCN42} is equivalent to $\Delta_{m3} = 0$. The ``if part'' follows.

Now we prove the ``only if'' part.

Since $\Gamma_4 = 0$ and  $\Gamma'_4 = 0$ are necessary conditions for the uniform order four we assume that $\Gamma_4 = 0$, $\Gamma'_4 = 0$ and $\Delta_{m3} = 0$ for some $m = 1,...,s^*$. Choose $\bar{m}$ such that
$\Delta_{\bar{m}3} \neq 0$  and $\Delta_{m3} = 0$ for $m = \bar{m}+1, ..., s^*$. Hence
$$
E'(\alpha h) = h^4 \sum_{m=1}^{\bar{m}}\left.\frac{\delta f}{\delta u}(\sigma+c^*_mh,y_{c^*_mh})\right|_{\theta}\left(\Delta_{m3}\left(\alpha, \frac{\theta}{h}\right)y^{(3)}(0)\right)+O(h^5),
$$
for $\alpha\in[0,1]$.

Let $\bar{\alpha} \in [0, 1]$ such that $\Delta_{\bar{m}3}(\bar{\alpha}, \cdot) \neq 0$. Since $\Delta_{\bar{m}3}(\bar{\alpha}, \cdot) \neq 0$ inside the interval $(-c^*_{\bar{m}} , -c^*_{\bar{m} -1})$
(where $c^*_{\bar{m} -1} = 0$ if $\bar{m} = 1$), there exists an interval $[\bar{\omega} - \varepsilon, \bar{\omega} + \varepsilon] \subseteq (-c^*_{\bar{m}} , -c^*_{\bar{m} -1})$, $\varepsilon > 0$, such that
$$
\Delta_{\bar{m}3}(\bar{\alpha}, \omega) = 0, \quad\omega \in [\bar{\omega} - \varepsilon, \bar{\omega} + \varepsilon].
$$
Set $A = -\bar{\omega}$.

Consider the scalar linear RFDE \eqref{eq:2order} defined by $r = Aa^2$ and
\begin{equation}
f(t, \phi) = \phi(-At^2) + g(t), \quad (t, \phi) \in (-a, a) \times \mathcal{C}
\label{eqN_proofOrd4}
\end{equation}
where $a > 1$ and $g : (-a, a) \to \mathbb{R}$ is such that $t \mapsto t^3$ is a solution of the RFDE on
$[-1 - r, a)$. Since $f$ is linear with respect to the second argument it has derivative given by
$$
\begin{aligned}
\frac{\delta f(t, \phi)}{\delta u}\psi &= \lim_{\rho\to0}\frac{f(t, \phi+\rho\psi) - f(t,\phi)}{\rho} \\
&= \lim_{\rho\to0}\frac{f(t, \phi)+\rho f(t,\psi) - f(t,\phi)}{\rho}=f(t, \psi)\\
&= \lim_{\rho\to0}\frac{\phi(-At^2)+g(t) +\rho\psi(-At^2) - \phi(-At^2)-g(t)}{\rho} = \psi(-At^2),
\end{aligned}
$$
for $(t, \phi) \in (-a,a)\times\mathcal{C}$, $\psi \in \mathcal{C}$.

Let
$$
\begin{aligned}
t_0 &= -1,\\
\phi_0(\theta) &= (-1 + \theta)^3,\quad \theta \in [-r, 0]. 
\end{aligned}
$$
Then
$$
u = u(t_0, \phi_0)(t) = t^3,\quad t \in [-1 - r, a).
$$
For $t \in [-1, a)$ and $h \in [0, \bar{h}(t, x_t))$ such that:
$$
|2\bar{\omega}c^*_mt + \bar{\omega}(c^*_m)^2t^2| \leq \varepsilon,\quad m = 1,..., \bar{m}
$$
and $h = t^2$, we have
$$
-\frac{A(t + c^*_mh)^2}{h} \in [\bar{\omega} - \varepsilon, \bar{\omega} + \varepsilon],\quad m = 1,..., \bar{m}
$$
and then the local error is given by (since $y^{(3)}(0)=6$)
$$
E'(\bar{\alpha}h) = 6h^4\Delta_{\bar{m}3}\left(\bar{\alpha}, -\frac{A(t +hc^*_{\bar{m}} h)^2}{h} \right) + O(h^5)
$$
with
$$
\Delta_{\bar{m}3} \left(\bar{\alpha}, -\frac{A(t +hc^*_{\bar{m}} h)^2}{h} \right) \geq \min_{\omega\in[\bar{\omega}-\varepsilon,\bar{\omega}+\varepsilon]} |\Delta_{\bar{m}3}(\bar{\alpha}, \omega)| > 0.
$$

Thus for the particular RFDE \eqref{eqN_proofOrd4}, for $t_0 = -1$, for $\phi_0(\theta) = (-1 + \theta)^3$, $\theta \in [-r, 0]$, and for $T = 1$ we have that: for all $H > 0$ and $C > 0$ there exist
$t \in [t_0, \bar{t}(t_0, \phi_0)) = [-1, a)$ and $h \in [0, \min\{H, \bar{h}(t, x_t)\})$ such that $t + \bar{c}h \leq T$ and
$$
\max_{\alpha\in[0,1]}|E'(\alpha h)| > Ch^5.
$$
So the method is not of uniform order 4.

\qed

Analogously for the discrete order four the following result holds.

\begin{theorem}
  An FCRKN method satisfying \eqref{OCN1} and \eqref{OCN2p} and  of discrete order three is of discrete order four iff
  \begin{equation}
   \sum_{i=1}^sb'_ic_i^3 = \frac{1}{4},
  \end{equation}
  \begin{equation}
   \sum_{\substack{i=1\\c_i=c^*_m}}^sb'_i\left(\sum_{j=1}^{s}a_{ij}(\beta)c_j - \frac{\beta^3}{6}\right)\!\!,\quad \beta\in[0,c^*_m]
  \end{equation}
  for $m=1,...,s^*$, and 
  \begin{equation}
   \sum_{i=1}^sb_ic_i^2 = \frac{1}{12}.
  \end{equation}
\end{theorem}

\subsection{Fifth order}
Now we develop the conditions for uniform and discrete orders four. As for order four we keep assuming that $f$ is of class $C^2$ with respect to second argument but now let us assume that $u$ is of piecewise class $C^6$.

Since \eqref{di_oh3} and \eqref{OCN2p}
$$
 E^i_{c_ih} = h^3\bar{\Gamma}_{i3}\left(\frac{\cdot}{h}\right)y^{(3)}(0) + h^4\bar{\Gamma}_{i4}\left(\frac{\cdot}{h}\right)y^{(4)}(0) + O(h^5),\quad i=1,...,s
$$
and then
\begin{equation}
 \begin{aligned}
  D_i &= h^3\left.\frac{\delta f}{\delta u}(\sigma+c_ih,y_{c_ih})\right|_{\theta}\left(\bar{\Gamma}_{i3}\left(\frac{\theta}{h}\right)y^{(3)}(0)\right) \\
  &\hspace{5mm} + h^4\left.\frac{\delta f}{\delta u}(\sigma+c_ih,y_{c_ih})\right|_{\theta}\left(\bar{\Gamma}_{i4}\left(\frac{\theta}{h}\right)y^{(4)}(0)\right) + O(h^5), \quad i=1,...,s.
 \end{aligned}
 \label{di_oh5}
\end{equation}

Now we have 
\begin{equation}
 \begin{aligned}
  E(\alpha h) &= h^3\Gamma_3(\alpha)y^{(3)}(0) + h^4\Gamma_4(\alpha)y^{(4)}(0) \\
  &\hspace{5mm}+ h^5 \sum_{m=1}^{s^*}\left.\frac{\delta f}{\delta u }(\sigma+c^*_mh,y_{c^*_mh})\right|_{\theta}\left(\sum_{\substack{i=1\\c_i=c^*_m}}^sb_i(\alpha)\bar{\Gamma}_{i3}\left(\frac{\theta}{h}\right)y^{(3)}(0)\right)\\
  &\hspace{5mm}+ h^5\Gamma_5(\alpha)y^{(5)}(0) + O(h^6), \quad \alpha\in[0,1],
 \end{aligned}
\end{equation}
and
\begin{equation}
 \begin{aligned}
  E'(\alpha h) &= h^2\Gamma'_2(\alpha)y^{(3)}(0) + h^3\Gamma'_3(\alpha)y^{(4)}(0) + h^4\Gamma'_4(\alpha)y^{(5)}(0) \\
  &\hspace{5mm}+ h^4 \sum_{m=1}^{s^*}\left.\frac{\delta f}{\delta u}(\sigma+c^*_mh,y_{c^*_mh})\right|_{\theta}\left(\sum_{\substack{i=1\\c_i=c^*_m}}^sb'_i(\alpha)\bar{\Gamma}_{i3}\left(\frac{\theta}{h}\right)y^{(3)}(0)\right)\\
  &\hspace{5mm}+ h^5 \sum_{m=1}^{s^*}\left.\frac{\delta f}{\delta u}(\sigma+c^*_mh,y_{c^*_mh})\right|_{\theta}\left(\sum_{\substack{i=1\\c_i=c^*_m}}^sb'_i(\alpha)\bar{\Gamma}_{i4}\left(\frac{\theta}{h}\right)y^{(4)}(0)\right)\\
  &\hspace{5mm}+ h^5\Gamma'_5(\alpha)y^{(6)}(0) + O(h^6), \quad \alpha\in[0,1].
 \end{aligned}
\end{equation}

The proof of the following theorems is analogous to the proof of Theorem\,\ref{TH_OCN4} with $g(t)$ taken such that a solution of the RFDE on $[-1-r, a)$ is $t \mapsto t^4$.

\begin{theorem}
  An FCRKN method satisfying \eqref{OCN1} and \eqref{OCN2p} and  of uniform order four is of uniform order five iff
  \begin{equation}
   \sum_{i=1}^sb'_i(\alpha)c_i^4 = \frac{\alpha^5}{5},\quad\alpha\in[0,1],
  \end{equation}
  \begin{equation}
   \sum_{\substack{i=1\\c_i=c^*_m}}^sb'_i(\alpha)\left(\sum_{j=1}^{s}a_{ij}(\beta)c_j^2 - \frac{\beta^4}{12}\right)\!\!,\quad\alpha\in[0,1],\quad \beta\in[0,c^*_m]
  \end{equation}
  for $m=1,...,s^*$,
  \begin{equation}
   \sum_{i=1}^sb_i(\alpha)c_i^3 = \frac{\alpha^5}{20},\quad\alpha\in[0,1],
  \end{equation}
  and 
  \begin{equation}
   \sum_{\substack{i=1\\c_i=c^*_m}}^sb_i(\alpha)\left(\sum_{j=1}^{s}a_{ij}(\beta)c_j - \frac{\beta^3}{6}\right)\!\!,\quad\alpha\in[0,1],\quad \beta\in[0,c^*_m]
  \end{equation}
  for $m=1,...,s^*$.
\end{theorem}

\begin{theorem}
  An FCRKN method satisfying \eqref{OCN1} and \eqref{OCN2p} and  of discrete order four is of discrete order five iff
  \begin{equation}
   \sum_{i=1}^sb'_ic_i^4 = \frac{1}{5},
  \end{equation}
  \begin{equation}
   \sum_{\substack{i=1\\c_i=c^*_m}}^sb'_i\left(\sum_{j=1}^{s}a_{ij}(\beta)c_j^2 - \frac{\beta^4}{12}\right)\!\!,\quad \beta\in[0,c^*_m],\quad m=1,...,s^*,
  \end{equation}
  \begin{equation}
   \sum_{i=1}^sb_ic_i^3 = \frac{1}{20},
  \end{equation}
  and 
  \begin{equation}
   \sum_{\substack{i=1\\c_i=c^*_m}}^sb_i\left(\sum_{j=1}^{s}a_{ij}(\beta)c_j - \frac{\beta^3}{6}\right)\!\!,\quad \beta\in[0,c^*_m], \quad m=1,...,s^*.
  \end{equation}
\end{theorem}

\section{FCRKNs with Reuse \label{sec:fcrkn_reuse}}
The methods constructed in~\cite{EJQTDE2015} satisfy the condition $b_i(\alpha) = \int b'_i(\alpha)d\alpha$, $i=1,...,s$. It isn't necessary but was used  because it instantly makes all order conditions for $b_i(\alpha)$ true, and could be satisfied with the minimal number of stages required to resolve the conditions for $b'_i(\alpha)$.

However in the case of \RKN{} methods we do not need to approximate $\dot{u}(t)$ values in order to compute the solution through the step. This means that it is sufficient to resolve discrete order conditions for $b'_i(\alpha)$ and only those uniform order conditions, which contain coefficients $b_i(\alpha)$. This (for methods of order 3 and higher) can be made with one stage fewer than full uniform order methods require.

And as for FCRKs with reuse we can compute the additional stage to be used for the uniform order approximation to $\dot{u}(t)$. This stage can be made equal to the first stage of the next step, if the following additional restrictions are satisfied:
\begin{itemize}    
  \item $c_s=1$;
  \item $b_i(\alpha) = a_{si}(\alpha)$ for any $i=1,...,s$;
  \item $\int_0^1 b'_i(\alpha)d\alpha = b_{i}(1)$ for any $i=1,...,s$.
\end{itemize}
The last condition provides continuous approximation to $\dot{u}$ over several steps.

FCRKNs with reuse are represented in Butcher tableaux of the form
\begin{equation}
  \begin{array}{c|ccccc}
    0    & &   \\
    c_2  & a_{21}(\alpha) \\
    \vdots & \vdots & \ddots  & \\
    c_{s-1}  & a_{s-1,1}(\alpha) & \cdots & a_{s-1,s-2}(\alpha) \\\hline
         & b_1(\alpha) & \cdots & b_{s-2}(\alpha) & b_{s-1}(\alpha) \\\hline
         & b'_1(\alpha) &  \cdots & b'_{s-2}(\alpha) & b'_{s-1}(\alpha) & b'_s(\alpha)
  \end{array}
  \label{tbl.frckn_gen}
\end{equation}

The following tables present methods of order 3 with 2 new stages per step and of order 4 with 4 new stages.

\subsection{Method of order 3}
\begin{equation}
  \begin{array}{c|ccc}
    0    & & &  \\
    \frac12    & \frac12\alpha^2\\\hline
    1    & \frac12\alpha^2-\frac13\alpha^3 & \frac13\alpha^3 \\\hline
           & \alpha-\frac32\alpha^2+\frac23\alpha^3 &  2\alpha^2-\frac43\alpha^3 & -\frac12\alpha^2+\frac23\alpha^3
  \end{array}
  \label{tbl.fcrkn3f}
\end{equation}

\begin{remark}
  Since the right-hand side of \eqref{eq:2order} doesn't depend on $\dot{u}$ a practically applicable method doesn't need a continuous extension for $\dot{u}$, i.e. it is sufficient to provide only discrete order $p$ for $\dot{u}$ with continuous order $p$ for $u$. For $p=3$ there exists the unique method with two stages (being an extension of the RKN method of order three for ODEs). However due to reuse its cost is the same as of the method above, which is more general and provides a continuous extension for $\dot{u}$ as well.
\end{remark}

\subsection{Method of order 4}
\begin{equation}
  \begin{array}{c|ccccc}
    0    & & & & & \\
    \frac{ 4}{11}  & \frac12\alpha^2\\
    \frac{10}{29}  & \frac12\alpha^2-\frac{11}{24}\alpha^3 & \frac{11}{24}\alpha^3 \\
    \frac{ 9}{11}  & \frac12\alpha^2-\frac{295}{696}\alpha^3& \frac{253}{232}\alpha^3& -\frac{2}{3}\alpha^3& & \\\hline
    1    & b_1(\alpha) & b_2(\alpha) & b_3(\alpha) & b_4(\alpha) \\\hline
         & b'_1(\alpha) & 0 & b'_3(\alpha) & b'_4(\alpha) & b'_5(\alpha)
  \end{array}
  \label{tbl.fcrkn4f}
\end{equation}
\begin{equation*}
  \begin{aligned}
    &b_1(\alpha) = \tfrac12\alpha^2 -\tfrac{5209361}{7811208}\alpha^3 +\tfrac{4299619}{15622416}\alpha^4
    &&b'_1(\alpha) = \alpha-\tfrac{461}{180}\alpha^2+\tfrac{23}{9}\alpha^3-\tfrac{319}{360}\alpha^4\\
    &b_2(\alpha) = \tfrac{960839}{1446520}\alpha^3-\tfrac{5770963}{8679120}\alpha^4
    &&b'_3(\alpha) = \tfrac{219501}{57380}\alpha^2-\tfrac{48778}{8607}\alpha^3+\tfrac{268279}{114760}\alpha^4\\
    &b_3(\alpha) = \tfrac{7}{43}\alpha^3+\tfrac{7}{43}\alpha^4
    &&b'_4(\alpha) = -\tfrac{6655}{2718}\alpha^2+\tfrac{17303}{2718}\alpha^3-\tfrac{38599}{10872}\alpha^4\\
    &b_4(\alpha) = -\tfrac{781726}{4882005}\alpha^3+\tfrac{4431163}{19528020}\alpha^4
    &&b'_5(\alpha) = \tfrac{45}{38}\alpha^2-\tfrac{371}{114}\alpha^3+\tfrac{319}{152}\alpha^4
  \end{aligned}
\end{equation*}

\section{Numerical Comparison}
To confirm the convergence order of the new methods we run multiple tests with constant step-size $h$ and measure the maximax error over the whole integration interval $Err = \max_{t_0\leq t \leq t_f}\|u(t)-\mathrm{H}(t)\|$, where $\mathrm{H}(t)$ is the continuous approximation to the solution by a numerical method. $Err$ should be proportional to $h^p$, where $p$ is the method's convergence order. We also compare number $N_f$ of right-hand sides $f$ evaluations required to provide certain $Err$. For FCRKNs we also measure $Err' = \max_{t_0\leq t \leq t_f}\|\dot{u}(t)-\mathrm{H}'(t)\|$, where $\mathrm{H}'(t)$ is the continuous approximation to the solution derivative.

\subsection{FCRKs}

We have chosen two DDE problems with overlapping for FCRKs to compare the methods \eqref{tbl.fcrk3f} and \eqref{tbl.fcrk4f} to the methods from~\cite{MasTorVer2005} of the same order.
  
\begin{figure}[!th]
\includegraphics[width=5.5cm]{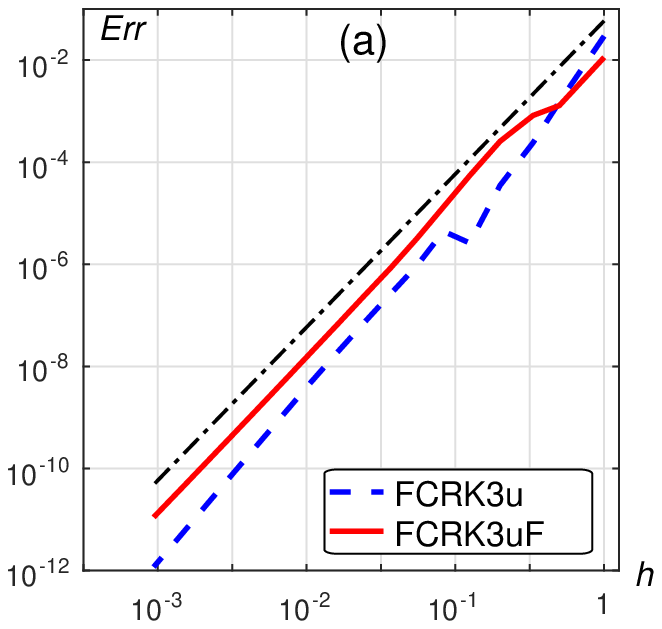}
\hspace{1cm}\includegraphics[width=5.5cm]{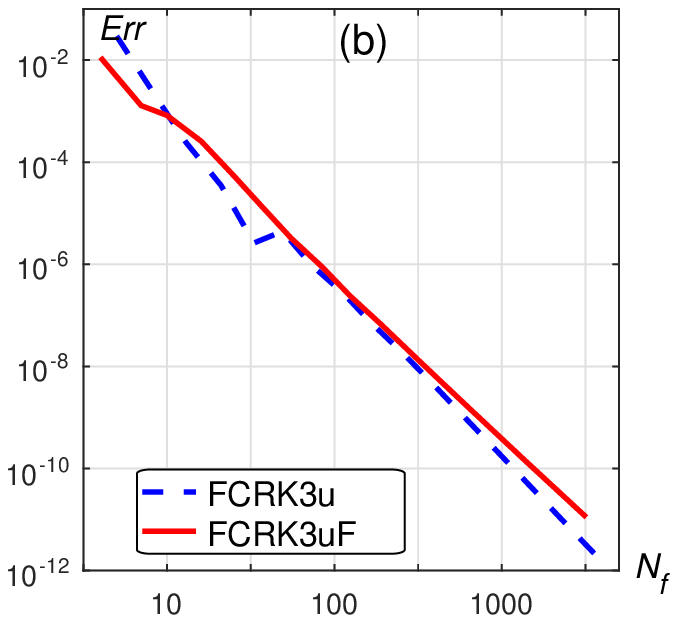}
\caption{Convergence orders (a) and error to $f$ evaluations (b) of order 3 methods for the Problem 1. The dot-dash reference line has slope 3. \label{fig_fcrk3p1}}
\end{figure}

\begin{figure}[!th]
\includegraphics[width=5.5cm]{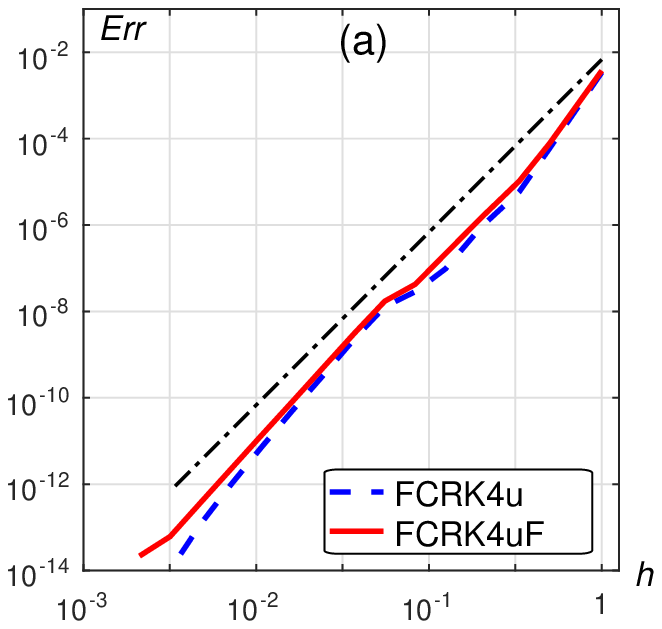}
\hspace{1cm}\includegraphics[width=5.5cm]{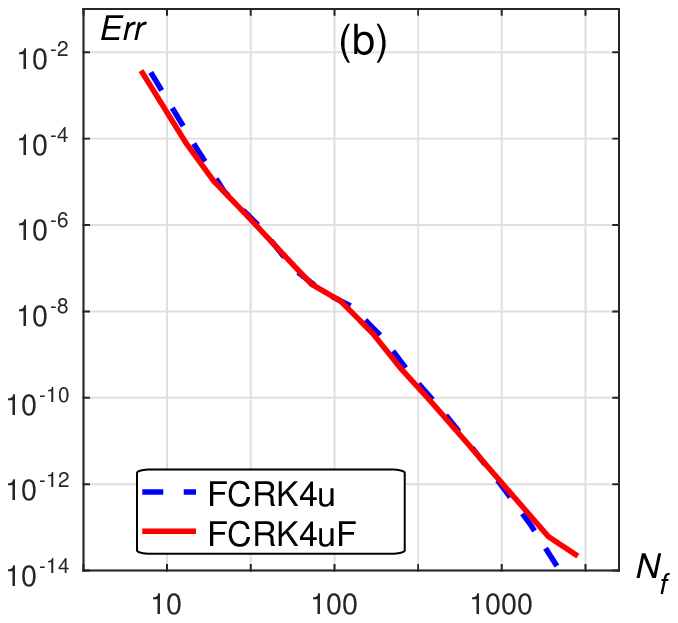}
\caption{Convergence orders (a) and error to $f$ evaluations (b) of order 4 methods for the Problem 1. The dot-dash reference line has slope 4. \label{fig_fcrk4p1}}
\end{figure}

\textbf{Problem 1} is the problem 1.2.6 from~\cite{PaulReport}. It is an initial value problem (IVP) with overlapping occuring for few first steps and no discontinuity  points:
\begin{equation}
  \left\{
  \begin{aligned}
    &\dot{u}(t) = u\left(\dfrac{t}{(1+2t)^2}\right)^{(1+2t)^2},\quad t\geq 0, \\
    &u(0)  = 1.  
  \end{aligned}\right.
  \label{Problem1}
\end{equation}
It has the analytical solution $u(t) = e^t$, $t\geq0$. We integrate \eqref{Problem1} at the interval $t\in[0, 1]$. The results are presented at Figs.~\ref{fig_fcrk3p1} and \ref{fig_fcrk4p1}.

\begin{figure}[!ht]
\includegraphics[width=5.5cm]{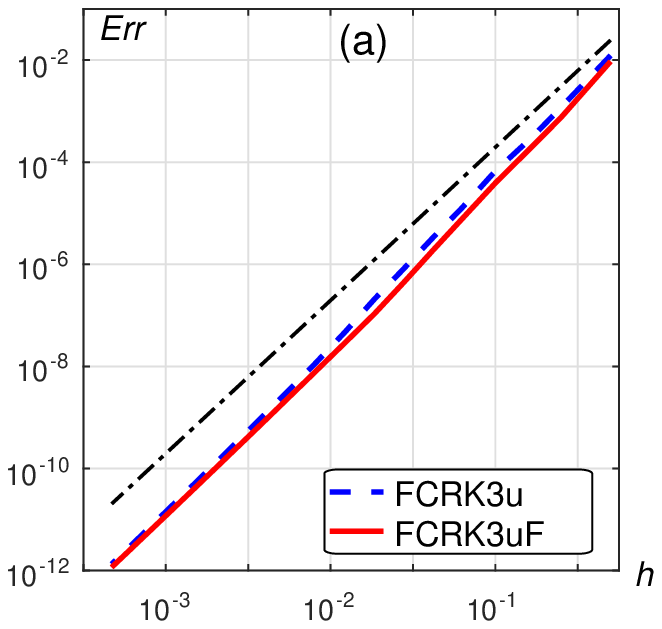}
\hspace{1cm}\includegraphics[width=5.5cm]{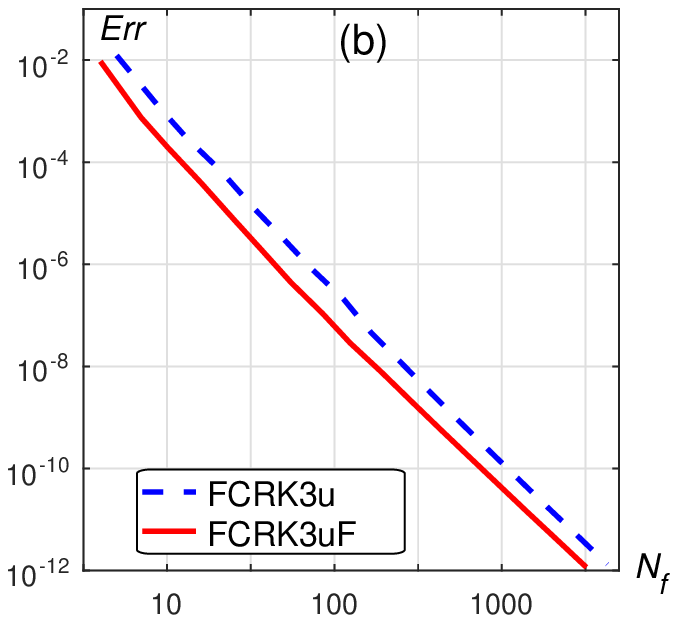}
\caption{Convergence orders (a) and error to $f$ evaluations (b) of order 3 methods for the Problem 2. The dot-dash reference line has slope 3. \label{fig_fcrk3p2}}
\end{figure}

\begin{figure}[!ht]
\includegraphics[width=5.5cm]{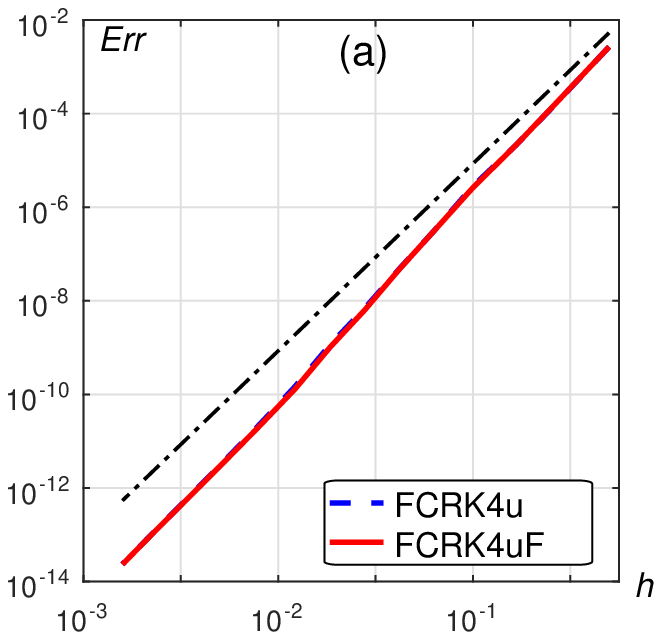}
\hspace{1cm}\includegraphics[width=5.5cm]{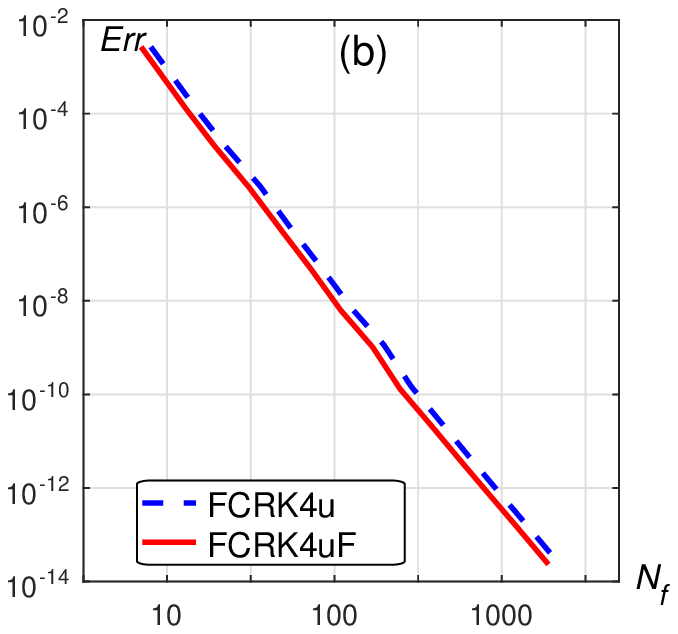}
\caption{Convergence orders (a) and error to $f$ evaluations (b) of order 4 methods for the Problem 2. The dot-dash reference line has slope 4. \label{fig_fcrk4p2}}
\end{figure}

\textbf{Problem 2} is the problem with vanishing delay. It has overlapping many times within the integration interval:
\begin{equation}
  \left\{
  \begin{aligned}
    &\dot{u}(t) = -u(g(t))u(t) e^{g(t)},&& t\geq 0, \\
    &u(t)  = e^{-t}, && t\leq 0.  
  \end{aligned}\right.
  \label{Problem2}
\end{equation}
where $g(t) = t-\frac{1}{100}\sin(100\pi t)^2$. It's analytical solution is the continuation of the history $u(t) = e^{-t}$. The problem is solved for $t\in[0, 0.5]$. The results are presented at Figs.~\ref{fig_fcrk3p2} and \ref{fig_fcrk4p2}.

As it can be seen both new methods show the expected convergence (at least for $h$ small enough). As for the computational costs, for methods with reuse they are in the most cases lower for the same global error.

\subsection{FCRKNs}

The problems of the second order are slightly modified problems \eqref{Problem1} and \eqref{Problem2}.
  
\begin{figure}[!ht]
\includegraphics[width=5.5cm]{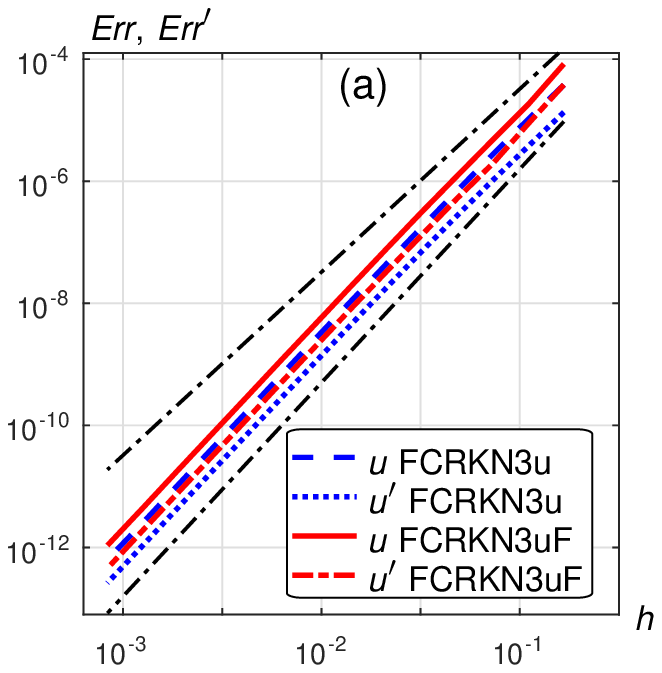}
\hspace{1cm}\includegraphics[width=5.5cm]{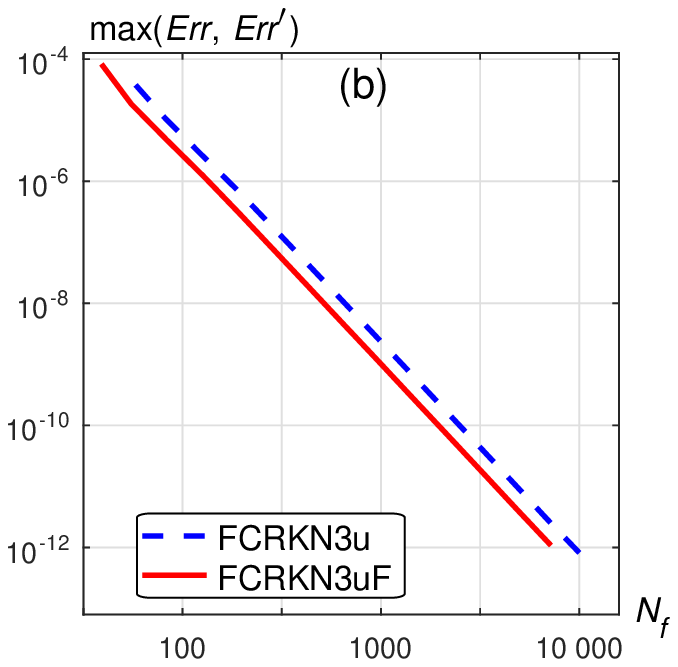}
\caption{Convergence orders (a) and maximal error to $f$ evaluations (b) of order 3 methods for the Problem 3. The dot-dash reference lines have slopes 3 and 3.5. 
\label{fig_fcrkn3p3}}
\end{figure}

\begin{figure}[!ht]
\includegraphics[width=5.5cm]{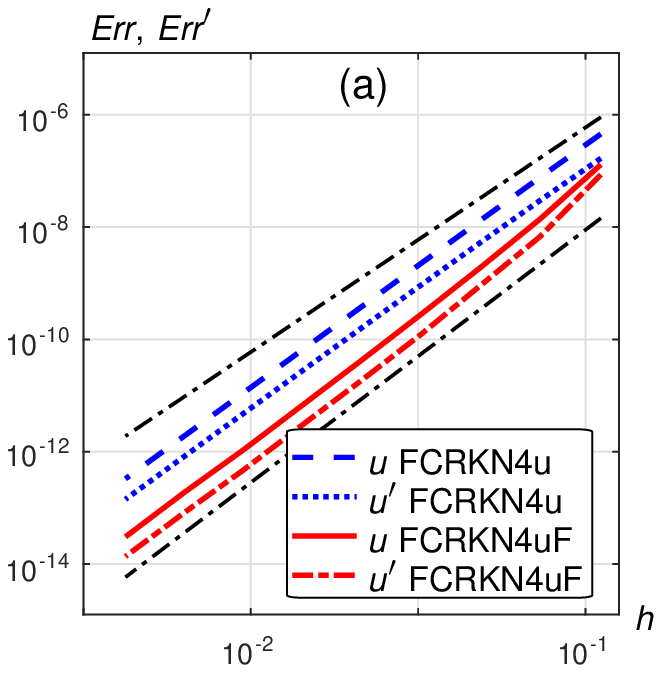}
\hspace{1cm}\includegraphics[width=5.5cm]{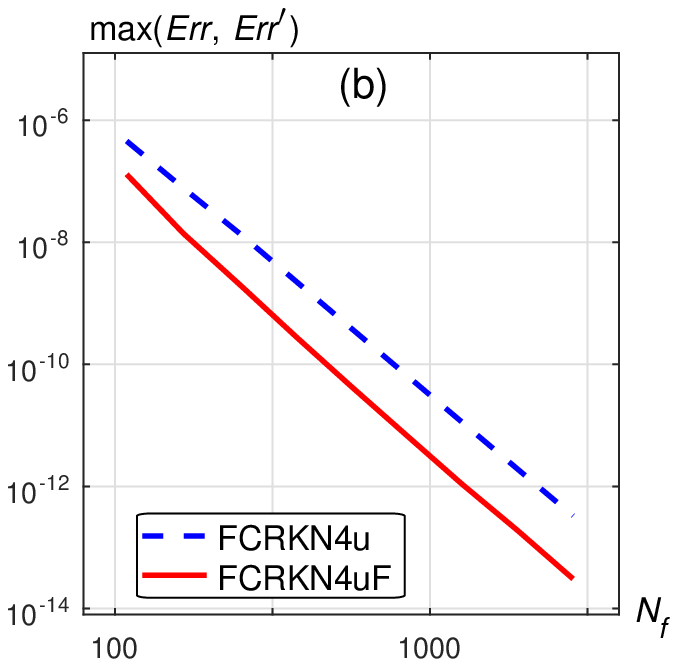}
\caption{Convergence orders (a) and maximal error to $f$ evaluations (b) of order 4 methods for the Problem 3. The dot-dash reference lines have slopes 4 and 4.5. \label{fig_fcrkn4p3}}
\end{figure}

\textbf{Problem 3} is an IVP as Problem 1:
\begin{equation}
  \left\{
  \begin{aligned}
    &\ddot{u}(t) = u\left(\dfrac{t}{(1+2t)^2}\right)^{(1+2t)^2},\quad t\geq 0, \\
    &u(0)  = 1, \quad \dot{u}(0) = -1
  \end{aligned}\right.
  \label{Problem3}
\end{equation}
has the solution $u(t) = e^{-t}$, $t\geq0$. We integrate it at the interval $t\in[0, 3]$. The results are presented at Figs.~\ref{fig_fcrkn3p3} and \ref{fig_fcrkn4p3}.

\begin{figure}[!ht]
\includegraphics[width=5.5cm]{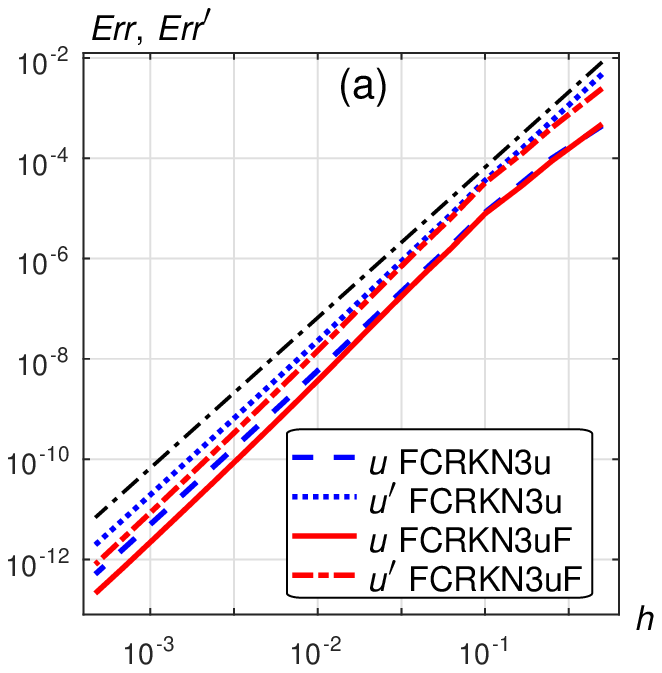}
\hspace{1cm}\includegraphics[width=5.5cm]{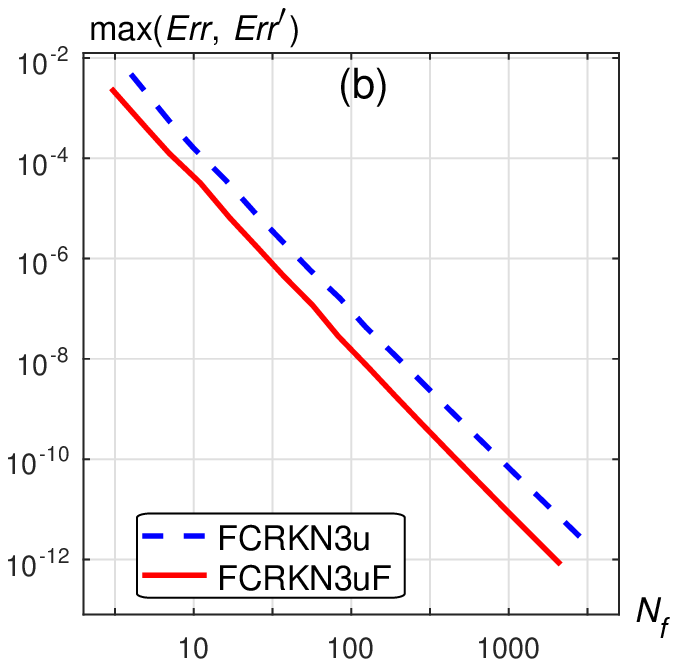}
\caption{Convergence orders (a) and maximal error to $f$ evaluations (b) of order 3 methods for the Problem 4. The dot-dash reference line has slope 3. 
\label{fig_fcrkn3p4}}
\end{figure}

\begin{figure}[!ht]
\includegraphics[width=5.5cm]{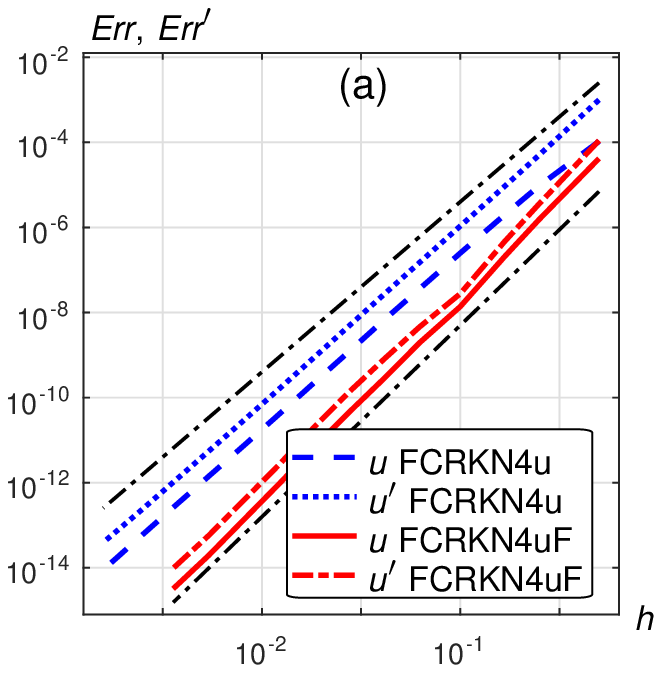}
\hspace{1cm}\includegraphics[width=5.5cm]{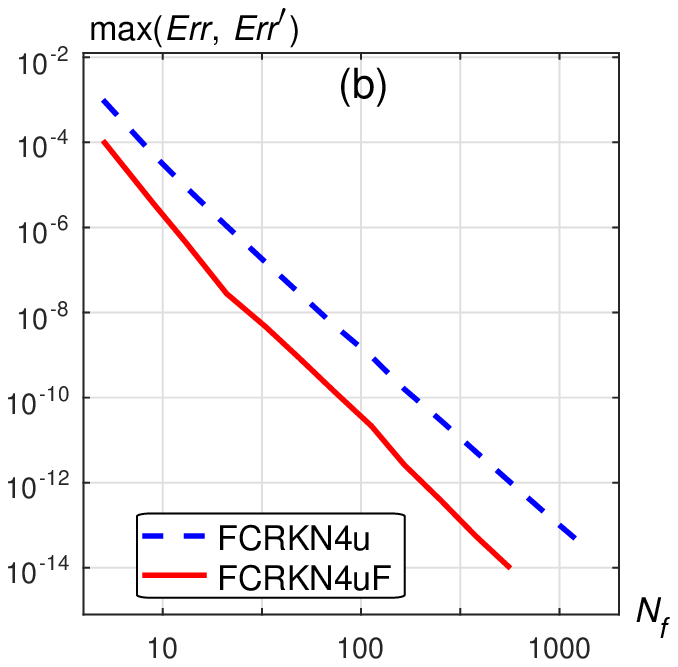}
\caption{Convergence orders (a) and maximal error to $f$ evaluations (b) of order 4 methods for the Problem 4. The dot-dash reference lines have slopes 4 and 4.5. \label{fig_fcrkn4p4}}
\end{figure}

\textbf{Problem 4} is based on Problem 2:
\begin{equation}
  \left\{
  \begin{aligned}
    &\ddot{u}(t) = u(g(t))u(t) e^{g(t)},&& t\geq 0, \\
    &u(t)  = e^{-t}, && t\leq 0.  
  \end{aligned}\right.
  \label{Problem4}
\end{equation}
with the same $g(t) = t-\frac{1}{100}\sin(100\pi t)^2$. The solution is $u(t) = e^{-t}$. The problem is solved for $t\in[0, 0.5]$. The results are presented at Figs.~\ref{fig_fcrkn3p4} and \ref{fig_fcrkn4p4}.

Here we also see that due to the reuse methods \eqref{tbl.fcrkn3f} and \eqref{tbl.fcrkn4f} require less computations. However, the convergence order is for some tests even higher than we've constructed, but less than a unit higher (about 3.5 for the third order methods and 4.5 for the method \eqref{tbl.fcrkn4f}). 

Similar results were observed in \cite{MagpantayPhD} and later explained in the unpublished talk \cite{HumphCortona12}. The matter is that if a method has order $p+1$ (or higher) without overlapping and order $p$ at the steps with overlapping its total convergence depends on the ratio of overlapping and no-overlapping steps. If the total length of the overlapping steps is $O(\sqrt{h})$ (which is the case in both problems) the convergence order is $p+0.5$.

\section{Conclusion}
Using the last stage of a Runge--Kutta type method as the first stage at the next step allows reducing the computational cost of continuous and functional continuous methods. We have considered first order RFDEs and second order RFDEs without dependency on the unknown function derivative. The constructed methods have the lowest possible number of stages for the uniform order they provide. The numerical 

\section{Acknowledgements}
The author would like to thank Prof. Stefano Maset for his valuable advices concerning functional continuous methods.

\section*{References}

\bibliographystyle{elsarticle-num}
\bibliography{eremin-anm-2018}

\end{document}